%% file: prime.tex
\documentclass{amsart}
\usepackage{amsmath,amssymb,amsthm}
\input rlepsf

\def\co{\colon\thinspace}
\newcommand{\Int}{\mbox{\rm Int}}

\newcommand{\R}{\mathbb{R}}

\newtheorem{thm}{Theorem}
\newtheorem{lem}[thm]{Lemma}

\theoremstyle{definition}
\newtheorem*{rem}{Remark}

\newtheorem*{ack}{Acknowledgements}

\begin{document}

\title[Unique decompostion for tight contact $3$-manifolds]{A unique
decomposition theorem for\\tight contact $3$-manifolds}
\author{Fan Ding}
\address{Department of Mathematics, Peking University,
Beijing 100871, P.~R. China}
\email{dingfan@math.pku.edu.cn}
\author{Hansj\"org Geiges}
\address{Mathematisches Institut, Universit\"at zu K\"oln,
Weyertal 86--90, 50931 K\"oln, F.R.G.}
\email{geiges@math.uni-koeln.de}
\date{}

\begin{abstract}
It has been shown by V.~Colin that every tight contact $3$-manifold
can be written as a connected sum of prime manifolds. Here we
prove that the summands in this decomposition are unique up to
order and contactomorphism.
\end{abstract}

\maketitle

\section{Introduction}
Unless stated otherwise, all $3$-manifolds in this note are assumed to
be closed, connected, and oriented.
A $3$-manifold is called {\bf non-trivial} if it is not
diffeomorphic to~$S^3$. A non-trivial $3$-manifold $P$ is said to be
{\bf prime} if in every connected sum
decomposition $P=P_0\# P_1$ one of the summands $P_0,P_1$
is~$S^3$. It is known that every non-trivial $3$-manifold
$M$ admits a prime decomposition, i.e.\ $M$ can be written as
a connected sum of finitely many prime manifolds. The main step in
the proof of this fact is due to H.~Kneser~\cite{knes29},
cf.~\cite{hemp76}. Moreover, as shown by J.~Milnor~\cite{miln62},
the summands in this prime decomposition are unique up to
order and diffeomorphism.

The purpose of the present note is to prove the analogous result
for tight contact $3$-manifolds. The basis for the argument
is a connected sum construction for such manifolds,
due to V.~Colin~\cite{coli97} and reproved by K.~Honda~\cite{hond02}.
Given a fixed connected sum decomposition $M=M_0\# M_1$ of
a $3$-manifold $M$, Colin's result says that
tight contact structures $\xi_i$ on $M_i$, $i=0,1$, give rise
to a tight contact structure $\xi_0\#\xi_1$ on $M$, uniquely defined
up to isotopy. Conversely, for any tight
contact structure $\xi$ on $M$ there are --- up to isotopy ---
unique tight contact structures $\xi_i$ on $M_i$, $i=0,1$,  such
that $\xi_0\#\xi_1$ is the given contact structure~$\xi$.
The prime decomposition theorem for tight contact $3$-manifolds is
an immediate consequence.

Although Colin's result goes a long way, it is not quite strong enough to
prove the unique decomposition
theorem for tight contact $3$-manifolds. This is due to the fact that
the system of $2$-spheres in a given manifold $M$ defining the prime
decomposition of $M$ is not, in general, unique
up to isotopy. The argument for the unique decomposition
of tight contact $3$-manifolds given here
closely follows the variant of Milnor's
argument given in J.~Hempel's book~\cite{hemp76}.
\section{Colin's results}
In this section we collect the results from \cite{coli97}
that we shall need. We assume that the reader is familiar
with the basics of contact topology at the level of
\cite{etny03}, \cite{hond04} and~\cite{geig06}.

\begin{lem}[{\cite[Lemme~5]{coli97}}]
\label{lem:0-surgery1}
Let $(M,\xi )$ be a (not necessarily connected)
tight contact $3$-manifold. Given embeddings
$f_0,f_1\co S^2\rightarrow M$, there is a contact structure
$\eta$ on $S^2\times [0,1]$ such that the characteristic foliation
$(S^2\times\{ i\})_{\eta}$ coincides with $S^2_{f_i^*\xi}$,
$i=0,1$.\footnote{Here $S^2_{f_i^*\xi}$ denotes
the characteristic foliation
induced by the embedding $f_i$, that is, the pull-back to $S^2$
via $f_i$ of the characteristic foliation $(f_i(S^2))_{\xi}$.}
This contact structure $\eta$ is unique up to isotopy rel boundary.\qed
\end{lem}

We can now define surgery along a $0$-sphere inside a given
(not necessarily connected) tight contact
$3$-manifold $(M,\xi )$ as follows; this includes the
formation of a connected sum.

Equip the $3$-disc $D^3$ with its standard orientation.
Let $\phi_0,\phi_1\co D^3\rightarrow M$ be embeddings such that
$\phi_0$ reverses and $\phi_1$ preserves orientation, and whose
images $B_i:=\phi_i(D^3)\subset M$ are disjoint.
Let $\eta$ be the contact structure on $S^2\times [0,1]$, constructed in
the preceding lemma, with the property that $(S^2\times\{ i\})_{\eta}=
(\partial D^3)_{\phi_i^*\xi}$. Then set
\[ (M',\xi')=(M\setminus\Int (B_0\cup B_1),\xi )\cup_{\partial}
(S^2\times [0,1],\eta ),\]
where $\cup_{\partial}$ denotes the obvious gluing along the boundary.

If $M=M_0+M_1$ is the topological sum of two connected tight
contact $3$-manifolds $(M_0,\xi_0)$, $(M_1,\xi_1)$, and
$B_i\subset M_i$, $i=0,1$, then $M'$ is the connected sum $M_0\#M_1$ of
$M_0$ and $M_1$, and we write $\xi_0\#\xi_1$ for the contact structure
$\xi'$ in this specific case. We also use the notation $(M_0,\xi_0)\# (M_1,
\xi_1)$ for this connected sum of tight contact $3$-manifolds.
As in the topological case, this connected sum operation is
commutative and associative;
these are consequences of the discussion that follows.
From Theorem~\ref{thm:Eliashberg} below it follows
that $(S^3,\xi_{st})$, the $3$-sphere with its unique tight
contact structure, serves as the neutral element.

\begin{lem}[{\cite[Corollaire~8]{coli97}}]
\label{lem:0-surgery2}
Let $(M',\xi')$ be a contact $3$-manifold and $f_t\co S^2\rightarrow M'$,
$t\in [0,1]$ an isotopy of embeddings. If the spheres $S_i:=f_i(S^2)$,
$i=0,1$, are convex with respect to $\xi'$, and $(M'\setminus S_0,\xi')$
is tight, then so is $(M'\setminus S_1,\xi')$.\qed
\end{lem}

\begin{lem}[{\cite[Proposition~9]{coli97}}]
\label{lem:0-surgery3}
The manifold $(M',\xi')$ obtained, in the way described above,
via $0$-surgery on a tight contact $3$-manifold $(M,\xi )$,
is tight and only depends,
up to contactomorphism, on the isotopy class
of the embeddings $\phi_0$, $\phi_1$.\qed
\end{lem}

\section{The unique decomposition theorem}
We can now formulate the unique decomposition theorem for tight
contact $3$-manifolds.

\begin{thm}
\label{thm:prime}
Every non-trivial tight contact $3$-manifold $(M,\xi )$ is
contactomorphic to a connected sum
\[ (M_1,\xi_1)\#\cdots\# (M_k,\xi_k)\]
of finitely many prime tight contact $3$-manifolds. The summands
$(M_i,\xi_i)$, $i=1,\ldots ,k$, are unique up to order and
contactomorphism.
\end{thm}

The proof of this theorem requires a few preparations.
Besides Colin's results, the most important ingredient is the
following theorem of Ya.~Eliashberg.

\begin{thm}[{\cite[Theorem~2.1.3]{elia92}}]
\label{thm:Eliashberg}
Two tight contact structures on the $3$-disc $D^3$ which induce the
same characteristic foliation on the boundary $\partial D^3$
are isotopic rel boundary.\qed
\end{thm}

First of all,
we observe that there is a well-defined
procedure for capping off a compact tight contact $3$-manifold with
boundary consisiting of a collection of $2$-spheres. Indeed, suppose
that $(M,\xi )$ is a tight contact $3$-manifold with boundary
$\partial M=S_1+\cdots +S_k$, where each $S_i$ is diffeomorphic to~$S^2$.
Choose orientation-reversing diffeomorphisms
$f_i\co\partial D^3\rightarrow S_i$.
By a reasoning as in Colin's proof of Lemma~\ref{lem:0-surgery1}, one finds
an orientation-preserving embedding $g_i\co D^3\rightarrow\R^3$
such that $S^2_{g_i^*\xi_{st}}=
S^2_{f_i^*\xi}$; here $\xi_{st}$ denotes the standard tight
contact structure on $\R^3$ (which is the restriction of $\xi_{st}$
on $S^3$ to the complement of a point).
The tight contact structure $\eta_i:=g_i^*\xi_{st}$ on $D^3$ ---
which by Theorem~\ref{thm:Eliashberg} is uniquely determined
by the characteristic foliation it induces on the boundary ---
can then be used to form the closed contact manifold
\[ (\widehat{M},\hat{\xi})=(M,\xi)\cup_{\partial}\bigl(
(D^3,\eta_1)\cup\ldots\cup (D^3,\eta_k)\bigr) ,\]
where the gluing is defined by the embeddings $f_i$.

Eliashberg's theorem entails that
we arrive at a contactomorphic
manifold if instead of gluing discs along the $S_i$ we first perturb
the boundary spheres into convex spheres $S_i'$ in the interior of $(M,\xi )$,
cut off the spherical shell between $S_i$ and $S_i'$, and then glue discs
along the $S_i'$. In the same way that Lemma~\ref{lem:0-surgery2}
enters the proof of Lemma~\ref{lem:0-surgery3} in Colin's argument,
one can use it here to conclude that $(\widehat{M},\hat{\xi})$ is tight.

\vspace{2mm}

Given an embedded $2$-sphere $S\subset\Int (M)$, we can find a product
neighbourhood $S\times [-1,1]\subset M$ of $S\equiv S\times\{ 0\}$.
Set $M_S=M\setminus \bigl( S\times (-1,1)\bigr)$. Again by
Theorem~\ref{thm:Eliashberg}, the contactomorphism type
of $(\widehat{M}_S,\hat{\xi})$ is independent of the choice of this
product neighbourhood; this follows by comparing the resulting
manifolds using two given product
neighbourhoods with a third manifold constructed from
a product neighbourhood contained in the first two. In particular,
this justifies our notation $(\widehat{M}_S,\hat{\xi})$.

\begin{lem}
\label{lem:prime1}
If $S_0$ and $S_1$ are isotopic $2$-spheres in $\Int (M)$, then
$(\widehat{M}_{S_0},\hat{\xi})$ and $(\widehat{M}_{S_1},\hat{\xi})$
are contactomorphic.
\end{lem}

\begin{proof}
This is clear if $S_1$ is isotopic to $S_0$ inside a product
neighbourhood $S_0\times (-1,1)$. The general case follows by an
argument very similar to Colin's proof of Lemma~\ref{lem:0-surgery2}.
\end{proof}

Given a connected sum decomposition $M=M_0\# M_1$ of a closed, connected
$3$-manifold with a tight contact structure~$\xi$, let $S\subset M$ be
an embedded sphere defining this connected sum, i.e.\
$\widehat{M}_S=M_0+M_1$. The described constructions imply that
\[ (M,\xi )= (M_0,\hat{\xi}|_{M_0})\# (M_1,\hat{\xi}|_{M_1}).\]
So the topological prime decomposition of $M$ also gives us
a decomposition of $(M,\xi )$ into prime tight contact $3$-manifolds.
The only remaining issue is the uniqueness of this decomposition
up to contactomorphism of the summands.

\vspace{2mm}

A $3$-manifold $M$ is said to be {\bf irreducible} if every
embedded $2$-sphere bounds a $3$-disc in $M$. Clearly, irreducible
$3$-manifolds (except~$S^3$) are prime.
There is but one orientable prime $3$-manifold
that is not irreducible, namely, $S^2\times S^1$ \cite[Lemma 3.13]{hemp76}.
In a connected sum $M=M_0\# S^2\times S^1$ we obviously find an embedded
non-separating $2$-sphere $S$ such that $\widehat{M}_S=M_0$;
simply take $S$ to be a fibre of $S^2\times S^1$ not affected by the
connected sum construction.

In the argument proving that the number of summands $S^2\times S^1$ in
a prime decomposition of $M$ is uniquely determined by $M$, the crucial
lemma is that for any two non-separating $2$-spheres $S_0,S_1\subset M$
there is a diffeomorphism of $M$ sending $S_0$ to $S_1$. In the
presence of a contact structure, this statement needs to be weakened slightly;
the following is sufficient for our purposes.

\begin{lem}
Let $(M,\xi )$ be a (connected) tight contact $3$-manifold
and $S_0,S_1\subset
M$ two non-separating $2$-spheres. Then
$(\widehat{M}_{S_0},\hat{\xi})$ and $(\widehat{M}_{S_1},\hat{\xi})$
are contactomorphic.
\end{lem}

\begin{proof}
By the preceding lemma we may assume that $S_0$ and $S_1$ are in
general position with respect to each other, so that $S_0\cap S_1$
consists of a finite number of embedded circles. We use induction on
the number $n$ of components of $S_0\cap S_1$.

If $n=0$, we may find disjoint product neighbourhoods
$S_i\times [-1,1]\subset M$, $i=0,1$. In case $M\setminus (S_0\cup S_1)$
is not connected, we may assume that the identifications
of these neighbourhoods with a product have
been chosen in such a way that $S_0\times\{ 1\}$ and $S_1\times\{ 1\}$
lie in the same component of $M\setminus (S_0\cup S_1)$.
As described above, we then obtain a well-defined
tight contact manifold $(\widetilde{M},\tilde{\xi})$
by capping off the boundary components $S_i\times\{\pm 1\}$ of
\[ M\setminus\bigl( S_0\times (-1,1)\cup S_1\times (-1,1)\bigr) \]
with $3$-discs $D_0^{\pm},D_1^{\pm}$. Our assumptions imply that
$D_0^-+D_0^+$ is isotopic to $D_1^-+D_1^+$ in $\widetilde{M}$.
By performing $0$-surgery with respect to these embeddings
of $S^0\times D^3$, we obtain $(\widehat{M}_{S_1},\hat{\xi})$
and $(\widehat{M}_{S_0},\hat{\xi})$, respectively, so the result
follows from Lemma~\ref{lem:0-surgery3}.

If $n>0$, then some component $J$ of $S_0\cap S_1$ bounds a
$2$-disc $D\subset S_1$ with $\Int (D)\cap S_0=\emptyset$.
Let $E'$ and $E''$ be the $2$-discs in $S_0$ bounded by $J$, and set
$S_0'=D\cup E'$ and $S_0''=D\cup E''$. At least one of $S_0'$ and $S_0''$,
say $S_0'$,
is non-separating.\footnote{Since $S_0$ is non-separating, there
is a loop $\gamma$ in $M$ (in general
position with respect to all spheres in question)
that intersects $S_0$ in a single point, say one
contained in the interior of $E'$. If $S_0''$ is separating,
then $\gamma$ intersects it in an even number of points. Since $\gamma$
does not intersect $E''$, these points all lie in $D$. So $\gamma$
intersects $S_0'$ in an odd number of points, which means that
$S_0'$ is non-separating.}
Move $S_0'$ slightly so that it becomes a smoothly embedded sphere
disjoint from $S_0$ and intersecting $S_1$ in fewer than $n$ circles.
Then two applications of the inductive assumption prove the inductive step.
\end{proof}

\begin{proof}[Proof of Theorem~\ref{thm:prime}]
As indicated above, it only remains to prove the uniqueness
statement. Thus, let
\[ (M_1,\xi_1)\# \cdots\# (M_k,\xi_k)\]
and
\[ (M_1^*,\xi_1^*)\# \cdots\# (M_l^*,\xi_l^*)\]
be two prime decompositions of a given tight contact $3$-manifold
$(M,\xi )$. Without loss of generality we
assume\footnote{Of course,
from the topological prime decomposition theorem, one already knows that
$k=l$, but this does not help to simplify the present proof.}
$k\leq l$ and use induction on~$k$. For $k=1$ there is nothing to prove.
Now assume $k>1$ and the assumption to be proved for
prime decompositions with fewer than $k$ summands.

\vspace{1mm}

(i) Suppose some $M_i$ (say $M_k$) is diffeomorphic
to $S^2\times S^1$. Then $M$ contains a non-separating $2$-sphere.
By applying the argument in the footnote to the preceding proof
to this non-separating $2$-sphere and the $2$-spheres defining
the splitting of $M$ into the connected sum of the $M_j^*$,
one finds a non-separating $2$-sphere in at least one of these
summands, say $M_l^*$, which therefore must be a copy of $S^2\times S^1$.
By a folklore theorem of Eliashberg,
there is a unique tight contact structure on~$S^2\times S^1$, cf.\
\cite{etny03} for an outline proof and \cite{geig} for a complete proof.
Thus, $(M_k,\xi_k)$ is contactomorphic to $(M_l^*,\xi_l^*)$.
Let $S_0,S_1$ be a fibre in $M_k, M_l^*$, respectively.
From Theorem~\ref{thm:Eliashberg} it follows that
\[ (\widehat{M}_{S_0},\hat{\xi})=(M_1,\xi_1)\#\cdots\# (M_{k-1},\xi_{k-1})\]
and
\[ (\widehat{M}_{S_1},\hat{\xi})=(M_1^*,\xi_1^*)\#
\cdots\# (M_{l-1}^*,\xi_{l-1}^*),\]
and by the preceding lemma these two manifolds are contactomorphic.
So the conclusion of the theorem follows from the inductive assumption.

\vspace{2mm}

(ii) It remains to deal with the case where all the $M_i$ are irreducible.
Arguing as before (with the roles of the two connected sum decompositions
reversed), we see that each $M_j^*$ must be irreducible. Choose
a separating $2$-sphere $S\subset M$ such that the closures $U,V$
of the components of $M\setminus S$ satisfy
\[ (\widehat{U},\widehat{\xi|_U})=(M_1,\xi_1)\#\cdots\# (M_{k-1},\xi_{k-1})\]
and $(\widehat{V},\widehat{\xi|_V})=(M_k,\xi_k)$.\footnote{The
contact structure $\widehat{\xi|_U}$ is the same as the restriction of
the contact structure $\hat{\xi}$ (defined on $\widehat{M}_S=
\widehat{U}+\widehat{V}$) to $\widehat{U}$.}

Similarly, there exist pairwise disjoint $2$-spheres $T_1,\ldots ,T_{l-1}$
in $M$ such that --- with $W_1,\ldots ,W_l$ denoting the closures of the
components of $M\setminus (T_1\cup\ldots\cup T_{l-1})$,
and $\xi_j$ the restriction of $\xi$ to $W_j$ --- we have
$(\widehat{W}_j,\hat{\xi}_j)=(M_j^*,\xi_j^*)$, $j=1,\ldots ,l$.

Suppose that the system $T_1,\ldots ,T_{l-1}$ of embedded spheres has been
chosen in general position with respect to $S$ and with
$S\cap (T_1\cup\ldots\cup T_{l-1})$ having the minimal number
of components among all such systems. 

Here we have to enter a caveat. The notation suggests that $W_1$ has
boundary $T_1$, the $W_j$ with $j\in\{ 2,\ldots ,l-1\}$ have boundary
$T_{j-1}\sqcup T_j$, and $W_l$ has boundary~$T_{l-1}$. In fact, some of
the reasoning in the proof given in \cite{hemp76} seems to rely on such
an assumption. However, under the minimality condition we have just
described, it is perfectly feasible that some of the $W_j$ have several
boundary components (i.e.\ the connected sum looks like a tree rather
than a chain). In particular, the numbering of the $W_j$ is not meant
to suggest any kind of order in which they are glued together.

We claim that the minimality condition implies
$S\cap (T_1\cup\ldots\cup T_{l-1})=\emptyset$.

Assuming this claim, we have $S\subset W_j$ for some
$j\in\{ 1,\ldots ,l\}$. Since $\widehat{W}_j=M_j^*$ is irreducible,
$S$ bounds a $3$-cell $B$ in $M_j^*$. Thus, $S$ cuts $W_j$ into
two pieces $X$ and $Y$, where say $\widehat{Y}=S^3$. By the
uniqueness of the tight contact structure on $S^3$ we have in fact
$(\widehat{Y},\widehat{\xi|_Y})=(S^3,\xi_{st})$. Moreover,
$(\widehat{X},\widehat{\xi|_X})=(M_j^*,\xi_j^*)$
by Theorem~\ref{thm:Eliashberg}.

Of the $3$-discs in $M_j^*$ used
for forming the connected sum with one or several of the other prime
manifolds, at least one has to be contained in~$B$, otherwise $S$ would
bound a disc in~$M$. This means that of the closures $U,V$
of the two components of $M\setminus S$, the one containing $Y$ must
contain at least one of $W_1,\ldots ,W_{j-1},W_{j+1},\ldots ,W_l$.
Thus, in the case $Y\subset V$, the numbering (including that
of $W_j$) can be chosen
in such a way that $W_1,\ldots ,W_{j-1},X\subset U$
and $Y,W_{j+1},\ldots ,W_l\subset V$, with $j\leq l-1$. (The case with
$X\subset V$ and $Y\subset U$ is analogous; here $j\geq 2$.)
With Theorem~\ref{thm:Eliashberg}, and in particular the fact that
$(S^3,\xi_{st})$ is the neutral element for the connected sum operation,
we conclude that
\begin{eqnarray*}
(M_1,\xi_1)\#\cdots\# (M_{k-1},\xi_{k-1}) & = & (\widehat{U},
                             \widehat{\xi|_U})\\
         & = & (\widehat{W}_1,\hat{\xi}_1)\#\cdots\#
                  (\widehat{W}_{j-1},\hat{\xi}_{j-1})
                 \# (\widehat{X},\widehat{\xi|X})\\
         & = & (M_1^*,\xi_1^*)\#\cdots \#(M_{j}^*,\xi_{j}^*)
\end{eqnarray*}
and
\begin{eqnarray*}
(M_k,\xi_k) & = & (\widehat{V},\widehat{\xi|_V})\\
            & = & (\widehat{Y},\widehat{\xi|_Y})\#
                   (\widehat{W}_{j+1},\hat{\xi}_{j+1})\#
                  \cdots\# (\widehat{W}_l,\hat{\xi}_l)\\
            & = & (M_{j+1}^*,\xi_{j+1}^*)\#\cdots\#
                  (M_l^*,\xi_l^*).
\end{eqnarray*}
Since $M_k$ is prime, we must have $j=l-1$, hence
$(M_k,\xi_k)=(M_l^*,\xi_l^*)$. Once again, the
theorem follows from the inductive assumption.

\begin{figure}[h]
\centerline{\relabelbox\small
\epsfxsize 10cm \epsfbox{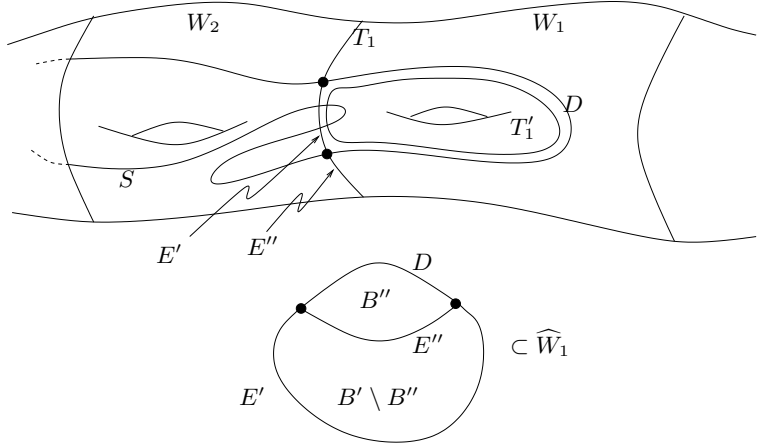}
\extralabel <-7.6cm,5.5cm> {$W_2$}
\extralabel <-5.4cm,5.3cm> {$T_1$}
\extralabel <-3.0cm,5.5cm> {$W_1$}
\extralabel <-8.5cm,3.4cm> {$S$}
\extralabel <-3.3cm,4.1cm> {$T_1'$}
\extralabel <-2.6cm,4.4cm> {$D$}
\extralabel <-4.6cm,2.3cm> {$D$}
\extralabel <-5.3cm,1.8cm> {$B''$}
\extralabel <-4.6cm,1.2cm> {$E''$}
\extralabel <-3.3cm,1.2cm> {$\subset\widehat{W}_1$}
\extralabel <-6.9cm,0.5cm> {$E'$}
\extralabel <-5.6cm,0.5cm> {$B'\setminus B''$}
\extralabel <-8.0cm,2.4cm> {$E'$}
\extralabel <-6.8cm,2.5cm> {$E''$}
\endrelabelbox}
\caption{Modification of the prime decomposition.}
\label{figure:prime}
\end{figure}

It remains to prove the claim. Arguing by contradiction, we assume
that $S\cap (T_1\cup\ldots\cup T_{l-1})\neq\emptyset$. Then we find
a $2$-disc $D\subset S$ with $\partial D\subset T_i$ for some
$i\in\{ 1,\ldots ,l-1\}$, and
$\Int (D)\cap (T_1\cup\ldots\cup T_{l-1})=\emptyset$.
This disc is contained in $W_j$ for some
$j\in\{ 1,\ldots ,l\}$. For ease of notation we assume that $i=j=1$,
and that $W_2$ is the other component adjacent to~$T_1$.

Let $E',E''$ be the $2$-discs in $T_1$ bounded by $\partial D$.
Since $\widehat{W}_1$ is irreducible, the sets
$D\cup E'$ and $D\cup E''$ (which are homeomorphic
copies of $S^2$) bound $3$-cells $B',B''$ in $\widehat{W}_1$.
One of these must contain the other, otherwise it would follow that
$\widehat{W}_1$ can be obtained by capping off the $3$-cell
$B'\cup_D B''$, and thus would be a $3$-sphere.

So suppose that $B''\subset B'$. Then $D\cup E'$ can be
deformed into a smooth $2$-sphere $T_1'$ that meets $S$ in fewer components
than $T_1$, see Figure~\ref{figure:prime}. In the complement
$M\setminus (T_1'\cup T_2\cup\ldots\cup T_{l-1})$ we still find
$W_3,\ldots ,W_l$, but $W_1,W_2$ have been changed to new components
$W_1'$, $W_2'$. Write $\xi_1',\xi_2'$, respectively, for the
restriction of $\xi$ to these components. We are done if we can show that
\[ (\widehat{W}_i',\hat{\xi}_i')=(\widehat{W}_i,\hat{\xi}_i),\;\; i=1,2,\]
because then the new system of spheres $T_1',T_2,\ldots ,T_{l-1}$
contradicts the minimality assumption on $T_1,T_2,\ldots ,T_{l-1}$.

The $2$-sphere $T_1'$ is isotopic to $T_1$ in $\widehat{W}_1$: simply move
$D\subset T_1'$ across the ball $B''$ to $E''$. But beware that $T_1'$ need
not be isotopic to $T_1$ in $W_1$ or $M$. However, $B''$ lies on the
same side of $T_1$ as $W_1$, so $T_1'$ is isotopic to $T_1$ in
\[ \widehat{W_1\cup W_2}=\widehat{W_1'\cup W_2'}.\]
Cutting this latter manifold open along $T_1$ and then capping off with
discs gives the disjoint union of $(\widehat{W}_1,\hat{\xi}_1)$
and $(\widehat{W}_2,\hat{\xi}_2)$; cutting it open along $T_1'$
and capping off yields the disjoint union of
$(\widehat{W}_1',\hat{\xi}_1')$ and $(\widehat{W}_2',\hat{\xi}_2')$.
From Lemma~\ref{lem:prime1} it follows that the results of either procedure
are contactomorphic.

This was the last point we had to show in order to
conclude the proof of the unique decomposition theorem.
\end{proof}

\begin{rem}
There is obviously no unique decomposition theorem for overtwisted
contact $3$-manifolds. For instance, start with a connected sum of
two distinct prime tight contact $3$-manifolds. Now perform a Lutz twist
in one or the other summand, preserving the topology
of the manifold and the homotopy class
of the contact structure as a $2$-plane field. By Eliashberg's
classification~\cite{elia89} of overtwisted contact structures,
the resulting manifolds are contactomorphic, and we obviously have
two distinct connected sum decompositions.
\end{rem}

\begin{ack}
This paper was written during a stay of F.~D. at the Mathematical
Institute of the Universit\"at zu K\"oln, supported by the DAAD grant no.\
A/06/27941. H.~G.\ is partially supported by the DFG grant no.\
GE 1245/1-2 within the framework of the Schwerpunktprogramm 1154 ``Globale
Differentialgeometrie''.
\end{ack}


\end{document}

%% file: rlepsf.tex
\input epsf
\def\relabelbox{%
  \hbox\bgroup%
}%
\def\endrelabelbox{%
\egroup%
}%
\def\relabel #1#2 {%
  \special{ps:/a {} def}%
  \smash{\rlap{#2}}%
}%
\def\adjustrelabel <#1,#2> #3#4 {%
  \special{ps:/a {} def}%
  \smash{\rlap{\kern #1 \raise #2\hbox{#4}}}%
}%
\def\extralabel <#1,#2> #3 {\smash{\rlap{\kern #1 \raise #2\hbox{#3}}}}%